\documentclass[10pt,journal,twocolumn]{IEEEtran}

\usepackage{amsmath,amssymb,threeparttable,placeins,makecell,stfloats,cite}
\usepackage[dvips]{graphicx,epstopdf}
\usepackage[usenames]{color}
\usepackage{amsfonts}
\usepackage{latexsym}
\usepackage{subfigure}
\usepackage[justification=centering]{caption}
\usepackage{floatrow}
\newtheorem{theorem}{{{\textit{Theorem}}}}

\newtheorem{lemma}{{{\textit{Lemma}}}}
\newtheorem{corollary}{{{{\textit{Corollary}}}}}

\newtheorem{definition}{{{\textit{Definition}}}}

\newtheorem{example}{{{\textit{Example}}}}

\begin{document}


\title{\textcolor{red}{A Direct Construction of Z-Complementary Pairs Using Generalized Boolean Functions}}
\author{Avik~Ranjan~Adhikary,~\IEEEmembership{Member,~IEEE},
		Palash~Sarkar,~
       Sudhan~Majhi,~\IEEEmembership{Senior Member,~IEEE}
\thanks{Avik Ranjan Adhikary and Palash Sarkar are with the Department of Mathematics, Indian Institute of Technology Patna, India, e-mail: {\tt avik.pma13@iitp.ac.in; palash.pma15@iitp.ac.in}. Sudhan Majhi is with Department of Electrical Engineering, Indian Institute of Technology Patna, India, e-mail: {\tt smajhi@iitp.ac.in}.}
}
 \maketitle

\begin{abstract}
 The zero correlation zone (ZCZ) ratio, i.e., the ratio of the width of the ZCZ and the length of the sequence plays a major role in reducing interference in an asynchronous environment of communication systems. However, to the best of the authors knowledge, the highest ZCZ ratio for even-length binary Z-complementary pairs (EB-ZCPs) which are directly constructed using generalized Boolean functions (GBFs), is $2/3$. In this letter, we present a direct construction of EB-ZCPs through GBFs, which can achieve a ZCZ ratio of $3/4$. In general, the constructed EB-ZCPs are of length $2^{m-1}+2$ ($m \in \mathbb{Z}^+$), having a ZCZ width of $2^{m-2}+2^{\pi(m-3)}+1$ where $\pi$ is a permutation over $m-2$ variables.

\end{abstract}

\section{INTRODUCTION}
The concept of ``complementary pair" was introduced by Golay in 1951 \cite{golay1}. Golay complementary pairs (GCPs) comprises of sequences whose aperiodic autocorrelation sums (AACSs) are zero at each out-of-phase time shift \cite{golay1}. GCPs have been found for numerous applications which include channel estimation \cite{spasojevic}, lowering the peak-to-mean envelope power ratio (PMEPR) \cite{popovic,davis,paterson}, RADAR waveform designs \cite{doppler_gcp}, etc. However, one of the main drawbacks of the GCPs is its limited availability. Binary GCPs exist only for lengths that are of the form $2^m10^n26^l$, where $m$, $n$, $l\in \mathbb{Z}^+$.

In search of binary sequences pairs having similar properties to that of GCPs, Fan \textit{et al.} proposed binary Z-complementary pairs (ZCPs), in 2007 \cite{fan}. ZCPs are sequence pairs, having zero AACSs at each out-of-phase time-shift within a zone around the zero shift position, commonly termed as the zero correlation zone (ZCZ) \cite{fan}. ZCPs are available for even-lengths as well as odd-lengths \cite{fan}. To know more about odd-length binary ZCPs (OB-ZCPs) and even-length binary ZCPs (EB-ZCPs) readers can go through \cite{fan,li2010,zilong_obzcp,Avik_iwsda,zilong_ebzcp,chen,Avik_ebzcp}. Along with GCPs, EB-ZCPs and OB-ZCPs can also be used as initial sequences to construct complementary sets \cite{Avik_cs}, complete complementary codes \cite{shibu1,shibu2} and Z-complementary code sets \cite{palash1,palash2}.

In \cite{fan}, Fan \textit{et al.} conjectured that for EB-ZCPs of length $N$, where $N$ is even ($N\neq 2^m10^n26^l$), the maximum ZCZ width $Z_{max}\leq N-2$. Working towards the solution of this open problem, Liu \textit{et al.} made a remarkable breakthrough in 2014 and proved that for a length $N$ EB-ZCP, the maximum ZCZ width ($Z_{max}$) that can be achieved is $N-2$ \cite{zilong_ebzcp}. The authors in \cite{zilong_ebzcp}, systematically designed EB-ZCPs of length $2^{m+1}+2^m$, which have a ZCZ width of $2^{m+1}$, by truncating certain binary GCPs of length $2^m$. The ZCZ ratio, i.e., the ratio of the width of the ZCZ and the sequence length, was measured to be $2/3$ for this construction. However, the problem of constructing length $N$ EB-ZCPs systematically, which can achieve a ZCZ width of $N-2$, is still an unsolved problem. Also, since the construction of EB-ZCPs requires GCPs as initial sequences, it was not a direct construction. 

Searching for more general construction of EB-ZCPs having larger ZCZ widths, Chen introduced a generalized Boolean function (GBF) based construction of EB-ZCPs, in 2017 \cite{chen}. This is the only direct construction of EB-ZCPs till date as it does not require any special sequences at the initial stage. Although the construction procedure was different from that of \cite{zilong_ebzcp}, however the ZCZ ratio of the resultant EB-ZCPs are still capped to $2/3$ \cite{chen}. For $n\leq m-2$, $n \in \mathbb{Z}^+$, the reported EB-ZCPs in \cite{chen} are of length $2^{m-1}+2^n$ and have the ZCZ width of $2^{m-2}+2^n$ \cite{chen}. 

Motivated by the works of \cite{zilong_ebzcp},\cite{chen}, to increase the ZCZ ratio of EB-ZCPs, recently we have proposed EB-ZCPs of length $2^{m+1}10^n26^l+2,~m\geq1,$ which have a ZCZ width of $3\times 2^{m-1}10^n26^l+1$ \cite{Avik_ebzcp}. The resultant sequences have the ZCZ ratio of $3/4$ \cite{Avik_ebzcp}. However, the construction was not direct as GCPs with certain intrinsic structural properties has been used at the initial stage. Then the insertion method have been applied to those GCPs to get the resultant EB-ZCPs \cite{Avik_ebzcp}. 


In search of direct construction, we propose the construction of EB-ZCPs through GBFs. Like the construction proposed by Chen in \cite{chen}, this construction also does not require any special sequences at the initial stage. However, the proposed EB-ZCPs have a wider ZCZ than the EB-ZCPs reported in \cite{chen}. To be specific, the ZCZ ratio of the EB-ZCPs resulted by our proposed construction is $3/4$ whereas in \cite{chen} the ZCZ ratio is $2/3$. The constructed EB-ZCPs are of length $2^{m-1}+2,~m\geq3$, having a ZCZ width of $2^{m-2}+2^{\pi(m-3)}+1$, where $\pi$ is a permutation over $m-2$ variables. When $\pi(m-3)=m-3$, the asymptotic ZCZ ratio becomes $3/4$. For EB-ZCPs having lengths of the form $2^{m-1}+2$, this ZCZ ratio is the maximum till date. It is quite impressive that proposed EB-ZCPs have exact AACS magnitude of $4$ at each time-shift outside the ZCZ where the AACS value is non-zero.

The remaining paper is organized as follows. We introduce GBFs and EB-ZCPs in Section II. In Section III, the proposed construction of EB-ZCPs is discussed. Finally, concluding remarks are addressed in Section IV. 

\section{NOTATIONS AND DEFINITIONS}
These notations will be followed throughout this paper. $\forall$ denotes `for all'. $+$ and $-$ denote $1$ and $-1$, respectively. Whenever it is not mentioned, binary sequences are sequences over $\mathbb{Z}_2$. We denote by $\bar{x}=1-x$
the binary complement of $x\in \{0,1\}$.

\vspace{0.1in}

\begin{definition}
Let ${\mathbf{a}}$ and ${\mathbf{b}}$ be two binary sequences of length $N$ over $\mathbb{Z}_2$. Then, the aperiodic cross-correlation function (ACCF) at a time-shift $\tau$ is defined by
\begin{equation}\label{defi_ACCF}
\rho_{\mathbf{a},\mathbf{b}}(\tau):= \left \{
\begin{array}{cl}
\sum\limits_{k=0}^{N-1-\tau}(-1)^{a_k+b_{k+\tau}},&~~0\leq \tau \leq N-1;\\
\sum\limits_{k=0}^{N-1-\tau}(-1)^{a_{k+\tau}+b_k},&~~-(N-1)\leq \tau \leq -1;\\
0,& ~~\mid \tau \mid \geq N.
\end{array}
\right .
\end{equation}
When the two sequences are identical, i.e., $\mathbf{a} = \mathbf{b}$, $\rho_{\mathbf{a},\mathbf{b}}(\tau)$ is known as an aperiodic auto-correlation function (AACF) of $\mathbf{a}$ and it is denoted by $\rho_{\mathbf{a}}(\tau)$.
\end{definition}

\vspace{0.1in}
\begin{definition}[EB-ZCP \cite{fan}, \cite{zilong_ebzcp}]
	Let $\mathbf{a}$ and $\mathbf{b}$ be two length $N$ binary sequences, where $N$ is even. ($\mathbf{a},\mathbf{b}$) is said to be an EB-ZCP with ZCZ width $Z$ iff
	\begin{equation}
	\rho_{\mathbf{a}}(\tau)+\rho_{\mathbf{b}}(\tau)=0,~~\forall~1\leq \tau \leq Z-1.
	\end{equation}
\end{definition}

\vspace{0.1in}
\begin{lemma}[\cite{zilong_ebzcp}]
	For a length $N$ EB-ZCP $(\mathbf{a},\mathbf{b})$, the maximum width of ZCZ, i.e., $Z_{max}\leq N-2$.
\end{lemma}

\vspace{0.1in}
\begin{lemma}[\cite{zilong_ebzcp}]
	Consider a length $N$ EB-ZCP $(\mathbf{a},\mathbf{b})$, which have a ZCZ width $Z<N$. Then 
	\begin{equation}
	\mid \rho_{\mathbf{a}}(\tau)+\rho_{\mathbf{b}}(\tau) \mid \geq 4, \forall~ Z \leq \tau < N,
	\end{equation}
	if it is non-zero.
\end{lemma}

\subsection{\textcolor{red}{Generalized Boolean Functions}}

A GBF $f:\mathbb{Z}_2^m\rightarrow\mathbb{Z}_{\textcolor{red}{q}}$ can uniquely be written as a linear combination of $2^m$ monomials
\begin{equation}
\begin{split}
1,x_0,x_1,\cdots,x_{m-1},&x_0x_1,x_0x_2,\cdots,\\&x_{m-2}x_{m-1}\cdots, x_0x_1\cdots x_{m-1},
\end{split}
\end{equation}
where the coefficients are taken from $\mathbb{Z}_{\textcolor{red}{q}}$.

By the notation $\Psi(f)$, we denote the sequence corresponding to a GBF $f$ and defined by $({\textcolor{red}{\omega}}^{f_0}, {\textcolor{red}{\omega}}^{f_1}, \cdots, {\textcolor{red}{\omega}}^{f_{2^m-1}})$ where \textcolor{red}{$\omega=\exp(2\pi\sqrt{-1}/q)$ ($q \geq 2$, is a positive integer),} $f_i=f(r_{i,0},r_{i,1},\dots, r_{i,m-1})$ and $\mathbf{r}_i\equiv(r_{i,0},r_{i,1},\dots, r_{i,m-1})$ is the binary vector representation of integer $i$ $(i=\displaystyle \sum_{j=0}^{m-1}r_{i,j}2^j)$. 


Given a GBF $f$ with $m$ variables, as defined above, the corresponding sequence $\Psi(f)$ will be of length $2^m$. In this paper, we concern about $(N,Z)$ ZCPs, where $N \neq 2^m$. Hence we define the truncated sequence $\Psi_{L}(f)$ corresponding to GBF $f$ by eliminating the first and last $L$ elements of the sequence $\Psi(f)$.

\begin{example}
	Let us consider $m=3$, $q=2$ and $f=x_0x_1+x_1x_2$, then
	\begin{equation}
		\begin{split}
		x_0x_1&=(0,0,0,1,0,0,0,1),\\
		x_1x_2&=(0,0,0,0,0,0,1,1),\\
		x_0x_1+x_1x_2&=(0,0,0,1,0,0,1,0),
		\end{split}
	\end{equation}
	and therefore $\Psi(f)=(+++-++-+)$. If we assume $L=1$, then $\Psi_{1}(f)=(++-++-)$. 
\end{example}

\begin{lemma}[\cite{rathinakumar}]\label{lem3}
	Consider a GBF $f:\mathbb{Z}_2^m\rightarrow \mathbb{Z}_{\textcolor{red}{q}}$, given by
	\begin{equation}
		f=\textcolor{red}{\frac{q}{2}}\sum_{\alpha=0}^{m-2}x_{\pi(\alpha)}x_{\pi(\alpha+1)}+\sum_{i=0}^{m-1}g_ix_i+g^\prime,
	\end{equation}
	and
	\begin{equation}
	\bar{f}=\textcolor{red}{\frac{q}{2}}\sum_{\alpha=0}^{m-2}\bar{x}_{\pi(\alpha)}\bar{x}_{\pi(\alpha+1)}+\sum_{i=0}^{m-1}g_i\bar{x}_i+g^\prime,
	\end{equation}
\end{lemma}
where $g_i,g^\prime \in \mathbb{Z}_q$. Then, $(\Psi(\bar{f}+x_{\pi(m-1)}),\Psi(\bar{f}))$ is one of the complementary mates of $(\Psi(f),\Psi(f+x_{\pi(m-1)}))$.

%
\section{Proposed \textcolor{red}{ZCPs} by using GBFs}

\subsection{Proposed Construction}

The proposed construction is discussed in this subsection.

\begin{theorem}\label{th1}
	For any integer $m\geq 4$, let $\pi$ be a permutation of $\{0,1,2,\dots, m-3\}$. For $d\in \mathbb{Z}_2$, let the GBF $g^d:\mathbb{Z}_2^m\rightarrow \mathbb{Z}_{\textcolor{red}{q}}$ be given as follows:
		\begin{equation}\label{gbf1}
		\begin{split}
g^d=\textcolor{red}{\frac{q}{2}} \textcolor{red}{[}x_{m-2}\bar{x}_{m-1}\cdot& \zeta^d+\bar{x}_{m-2}x_{m-1}\cdot \eta^d+d\bar{x}_{m-2}\bar{x}_{m-1}  \\ & 
+x_{m-2}x_{m-1}\textcolor{red}{]}\textcolor{red}{+\sum_{i=0}^{m-3}e_ix_i+\sum_{i=0}^{m-3}f_i\bar{x}_i}
		\end{split}
		\end{equation}
	where $\zeta^d:\mathbb{Z}_2^{m-3}\rightarrow \mathbb{Z}_2$ is 
	\begin{equation}
	\zeta^d=\sum_{\alpha=0}^{m-4}x_{\pi(\alpha)}x_{\pi(\alpha+1)}+dx_{\pi(m-3)},
	\end{equation}
	and $\eta^d:\mathbb{Z}_2^{m-3}\rightarrow \mathbb{Z}_2$ is
	\begin{equation}
	\eta^d=\sum_{\alpha=0}^{m-4}\bar{x}_{\pi(\alpha)}\bar{x}_{\pi(\alpha+1)}+\bar{d}x_{\pi(m-3)}.
	\end{equation}
	Then $(\mathbf{a},\mathbf{b})\equiv \left(\Psi_{2^{m-2}-1}(g^0),\Psi_{2^{m-2}-1}(g^1)\right)$ forms an EB-ZCP of length $2^m-2\times (2^{m-2}-1)=2^{m-1}+2$, having ZCZ width of $2^{m-2}+2^{\pi(m-3)}+1$. 
\end{theorem}

\begin{IEEEproof} 
		 For $0<\tau\leq2^{m-2}$, using \textit{Lemma} \ref{lem3}, we have
		
		\begin{equation}\label{eq11}
		\begin{split}
		&\rho_{\mathbf{a}}(\tau)+\rho_{\mathbf{b}}(\tau)\\&=\left[{\textcolor{red}{\omega}}^{\zeta^0(\mathbf{r}_{2^{m-2}+\tau-1})}-{\textcolor{red}{\omega}}^{\eta^0(\mathbf{r}_{3\times2^{m-2}-\tau})}\right]\\&\quad +\left[-{\textcolor{red}{\omega}}^{\zeta^1(\mathbf{r}_{2^{m-2}+\tau-1})}-{\textcolor{red}{\omega}}^{\eta^1(\mathbf{r}_{3\times2^{m-2}-\tau})}\right]\\
		&=\left[{\textcolor{red}{\omega}}^{\zeta^0(\mathbf{r}_{2^{m-2}+\tau-1})}-{\textcolor{red}{\omega}}^{\zeta^1(\mathbf{r}_{2^{m-2}+\tau-1})}\right]\\&\quad -\left[{\textcolor{red}{\omega}}^{\eta^0(\mathbf{r}_{3\times2^{m-2}-\tau})}+{\textcolor{red}{\omega}}^{\eta^1(\mathbf{r}_{3\times2^{m-2}-\tau})}\right].
		\end{split}
		\end{equation}
		Since,
		\begin{equation}\label{eq12}
		\zeta^1(\mathbf{r}_{2^{m-2}+\tau-1})=\zeta^0(\mathbf{r}_{2^{m-2}+\tau-1})+r_{2^{m-2}+\tau-1,\pi(m-3)},
		\end{equation}
		and
		\begin{equation}\label{eq13}
		\eta^0(\mathbf{r}_{3\times2^{m-2}-\tau})=\eta^1(\mathbf{r}_{3\times2^{m-2}-\tau})+r_{3\times2^{m-2}-\tau,\pi(m-3)},
		\end{equation}
		applying (\ref{eq12}) and (\ref{eq13}) in (\ref{eq11}), we get,
		
		\begin{equation}\label{eq14}
		\begin{split}
		&\rho_{\mathbf{a}}(\tau)+\rho_{\mathbf{b}}(\tau)\\&={\textcolor{red}{\omega}}^{\zeta^0(\mathbf{r}_{2^{m-2}+\tau-1})}\left[1-{\textcolor{red}{\omega}}^{r_{2^{m-2}+\tau-1,\pi(m-3)}}\right]\\& \quad -{\textcolor{red}{\omega}}^{\eta^1(\mathbf{r}_{3\times2^{m-2}-\tau})}\left[{\textcolor{red}{\omega}}^{r_{3\times2^{m-2}-\tau,\pi(m-3)}}+1\right].
		\end{split}
		\end{equation}
	So, we have the following two sub-cases:
		\begin{enumerate}
			\item For $r_{2^{m-2}+\tau-1,\pi(m-3)}=0$, we have $r_{3\times2^{m-2}-\tau,\pi(m-3)}=1$. In this case, we can easily conclude that 
			\begin{equation}
			\rho_{\mathbf{a}}(\tau)+\rho_{\mathbf{b}}(\tau)=0.
			\end{equation}
			
			\item For $r_{2^{m-2}+\tau-1,\pi(m-3)}=1$, we have $r_{3\times2^{m-2}-\tau,\pi(m-3)}=0$. In this case, we have from (\ref{eq14})
			
			\begin{equation}\label{eq16}
			\begin{split}
			&\rho_{\mathbf{a}}(\tau)+\rho_{\mathbf{b}}(\tau)\\&=2\left[{\textcolor{red}{\omega}}^{\zeta^0(\mathbf{r}_{2^{m-2}+\tau-1})}- {\textcolor{red}{\omega}}^{\eta^1(\mathbf{r}_{3\times2^{m-2}-\tau})}\right].
			\end{split}
			\end{equation}
			
			Since, $\Psi(\eta^1)$ is the reverse sequence of $\Psi(\zeta^0)$, therefore
			\begin{equation}
			{\textcolor{red}{\omega}}^{\zeta^0(\mathbf{r}_{2^{m-2}+\tau-1})}={\textcolor{red}{\omega}}^{\eta^1(\mathbf{r}_{3\times2^{m-2}-\tau})}.
			\end{equation}
			And hence
			\begin{equation}
			\rho_{\mathbf{a}}(\tau)+\rho_{\mathbf{b}}(\tau)=0.
			\end{equation}
		\end{enumerate}
		
		 For $2^{m-2} < \tau \leq 2^{m-2}+2^{\pi(m-3)}$, using \textit{Lemma} \ref{lem3}, we have
		
		\begin{equation}\label{eq19}
		\begin{split}
		&\rho_{\mathbf{a}}(\tau)+\rho_{\mathbf{b}}(\tau)\\&=\left[{\textcolor{red}{\omega}}^{\eta^0(\mathbf{r}_{2^{m-2}+\tau-1})}-{\textcolor{red}{\omega}}^{\zeta^0(\mathbf{r}_{3\times2^{m-2}-\tau})}\right]\\& \quad +\left[-{\textcolor{red}{\omega}}^{\eta^1(\mathbf{r}_{2^{m-2}+\tau-1})}-{\textcolor{red}{\omega}}^{\zeta^1(\mathbf{r}_{3\times2^{m-2}-\tau})}\right]\\
		&=\left[{\textcolor{red}{\omega}}^{\eta^0(\mathbf{r}_{2^{m-2}+\tau-1})}-{\textcolor{red}{\omega}}^{\eta^1(\mathbf{r}_{2^{m-2}+\tau-1})}\right]\\& \quad -\left[{\textcolor{red}{\omega}}^{\zeta^0(\mathbf{r}_{3\times2^{m-2}-\tau})}+{\textcolor{red}{\omega}}^{\zeta^1(\mathbf{r}_{3\times2^{m-2}-\tau})}\right].
		\end{split}
		\end{equation}
		Since,
		\begin{equation}\label{eq20}
		\zeta^1(\mathbf{r}_{3\times2^{m-2}-\tau})=\zeta^0(\mathbf{r}_{3\times2^{m-2}-\tau})+r_{3\times2^{m-2}-\tau,\pi(m-3)},
		\end{equation}
		and
		\begin{equation}\label{eq21}
		\eta^0(\mathbf{r}_{2^{m-2}+\tau-1})=\eta^1(\mathbf{r}_{2^{m-2}+\tau-1})+r_{2^{m-2}+\tau-1,\pi(m-3)},
		\end{equation}
		applying (\ref{eq20}) and (\ref{eq21}) in (\ref{eq19}), we get,
		
		\begin{equation}\label{eq22}
		\begin{split}
		&\rho_{\mathbf{a}}(\tau)+\rho_{\mathbf{b}}(\tau)\\&={\textcolor{red}{\omega}}^{\eta^1(\mathbf{r}_{2^{m-2}+\tau-1})}\left[{\textcolor{red}{\omega}}^{r_{2^{m-2}+\tau-1,\pi(m-3)}}-1\right]\\& \quad -{\textcolor{red}{\omega}}^{\zeta^0(\mathbf{r}_{3\times2^{m-2}-\tau})}\left[1+{\textcolor{red}{\omega}}^{r_{3\times2^{m-2}-\tau,\pi(m-3)}}\right].
		\end{split}
		\end{equation}
		Again, we have the following sub-case:
		\begin{enumerate}
			\item When $r_{2^{m-2}+\tau-1,\pi(m-3)}=0$, then $r_{3\times2^{m-2}-\tau,\pi(m-3)}=1$. In this case, we can easily conclude that 
			\begin{equation}
			\rho_{\mathbf{a}}(\tau)+\rho_{\mathbf{b}}(\tau)=0.
			\end{equation}	
		\end{enumerate}
		
		So, the ZCZ width is $2^{m-2}+2^{\pi(m-3)}+1$.
		

\end{IEEEproof}

\subsection{The AACS magnitude, outside the ZCZ}

\begin{corollary}
The obtained non-zero magnitude of AACS of the proposed EB-ZCPs in \textit{Theorem} \ref{th1} is exactly $4$ outside the ZCZ, if it is not zero.
\end{corollary}

\begin{IEEEproof}
When $\tau >2^{m-2}+2^{\pi(m-3)}$, we have the following cases:

\begin{enumerate}
	\item For $r_{2^{m-2}+\tau-1,\pi(m-3)}=0$, we have $r_{3\times2^{m-2}-\tau,\pi(m-3)}=1$. In this case, we can easily conclude that 
	\begin{equation}
	\rho_{\mathbf{a}}(\tau)+\rho_{\mathbf{b}}(\tau)=0.
	\end{equation}

	\item For $r_{2^{m-2}+\tau-1,\pi(m-3)}=1$, we have $r_{3\times2^{m-2}-\tau,\pi(m-3)}=0$. In this case, using (\ref{eq20}) and (\ref{eq21}), we have from (\ref{eq19})
	
	\begin{equation}\label{eq25}
	\begin{split}
	&\rho_{\mathbf{a}}(\tau)+\rho_{\mathbf{b}}(\tau)\\&={\textcolor{red}{\omega}}^{\eta^0(\mathbf{r}_{2^{m-2}+\tau-1})}\left[1-{\textcolor{red}{\omega}}^{r_{2^{m-2}+\tau-1,\pi(m-3)}}\right]\\& \quad -{\textcolor{red}{\omega}}^{\zeta^0(\mathbf{r}_{3\times2^{m-2}-\tau})}\left[1+{\textcolor{red}{\omega}}^{r_{3\times2^{m-2}-\tau,\pi(m-3)}}\right]\\
	&=2\left[{\textcolor{red}{\omega}}^{\eta^0(\mathbf{r}_{2^{m-2}+\tau-1})}-{\textcolor{red}{\omega}}^{\zeta^0(\mathbf{r}_{3\times2^{m-2}-\tau})}\right].
	\end{split}
	\end{equation}
	
	Since,
	\begin{equation}
	\eta^0(\mathbf{r}_{2^{m-2}+\tau-1})=\zeta^0(\mathbf{r}_{3\times2^{m-2}-\tau})+1,
	\end{equation}
	from (\ref{eq25}) we get 
	\begin{equation}\label{eq27}
	\begin{split}
	&\rho_{\mathbf{a}}(\tau)+\rho_{\mathbf{b}}(\tau)\\
	&=2\left[{\textcolor{red}{\omega}}^{\eta^0(\mathbf{r}_{2^{m-2}+\tau-1})}-{\textcolor{red}{\omega}}^{\zeta^0(\mathbf{r}_{3\times2^{m-2}-\tau})}\right]\\
	&=-4 \times {\textcolor{red}{\omega}}^{\zeta^0(\mathbf{r}_{3\times2^{m-2}-\tau})}.
	\end{split}
	\end{equation}
\end{enumerate}
Thus, $|\rho_{\mathbf{a}}(\tau)+\rho_{\mathbf{b}}(\tau)|=4$, which is minimum.
\end{IEEEproof}

\subsection{The ZCZ ratio}

When $\pi(m-3)=m-3$, then we get the maximum ZCZ width. And in that case

\begin{equation}
\begin{split}
	\frac{Z}{N}&=\lim\limits_{m\rightarrow\infty}\frac{2^{m-2}+2^{\pi(m-3)}+1}{2^{m-1}+2}\\
	&=\lim\limits_{m\rightarrow\infty}\frac{2^{m-2}+2^{m-3}+1}{2^{m-1}+2}\\&=\lim\limits_{m\rightarrow\infty}\frac{3}{4}-\frac{1}{2^m+4}\\
	&\approx \frac{3}{4}.
\end{split}
\end{equation}

\begin{example}\label{ex2}
	Let us consider $m=6$ and $\{\pi(0),\pi(1),\pi(2),\pi(3)\}=\{2,0,1,3\}$. For $d\in \mathbb{Z}_2$, let the Boolean function $g^d$ be given by 
	\begin{equation}
	g^d=x_4\bar{x}_5\cdot \zeta^d+\bar{x}_4x_5\cdot \eta^d+d\bar{x}_4\bar{x}_5+x_4x_5,
	\end{equation}
	where
	\begin{equation}
	\zeta^d=x_2x_0+x_0x_1+x_1x_3+dx_3,
	\end{equation}
	and
	\begin{equation}
	\eta^d=\bar{x}_2\bar{x}_0+\bar{x}_0\bar{x}_1+\bar{x}_1\bar{x}_3+\bar{d}x_3.
	\end{equation}
	Since $\pi(3)=3$, the pair $(\mathbf{a},\mathbf{b})\equiv \left(\Psi_{15}(g^0),\Psi_{15}(g^1)\right)$ gives a $34$ length EB-ZCP, with ZCZ width $25$.
\end{example}

\begin{figure}
	\includegraphics[width=0.9\columnwidth]{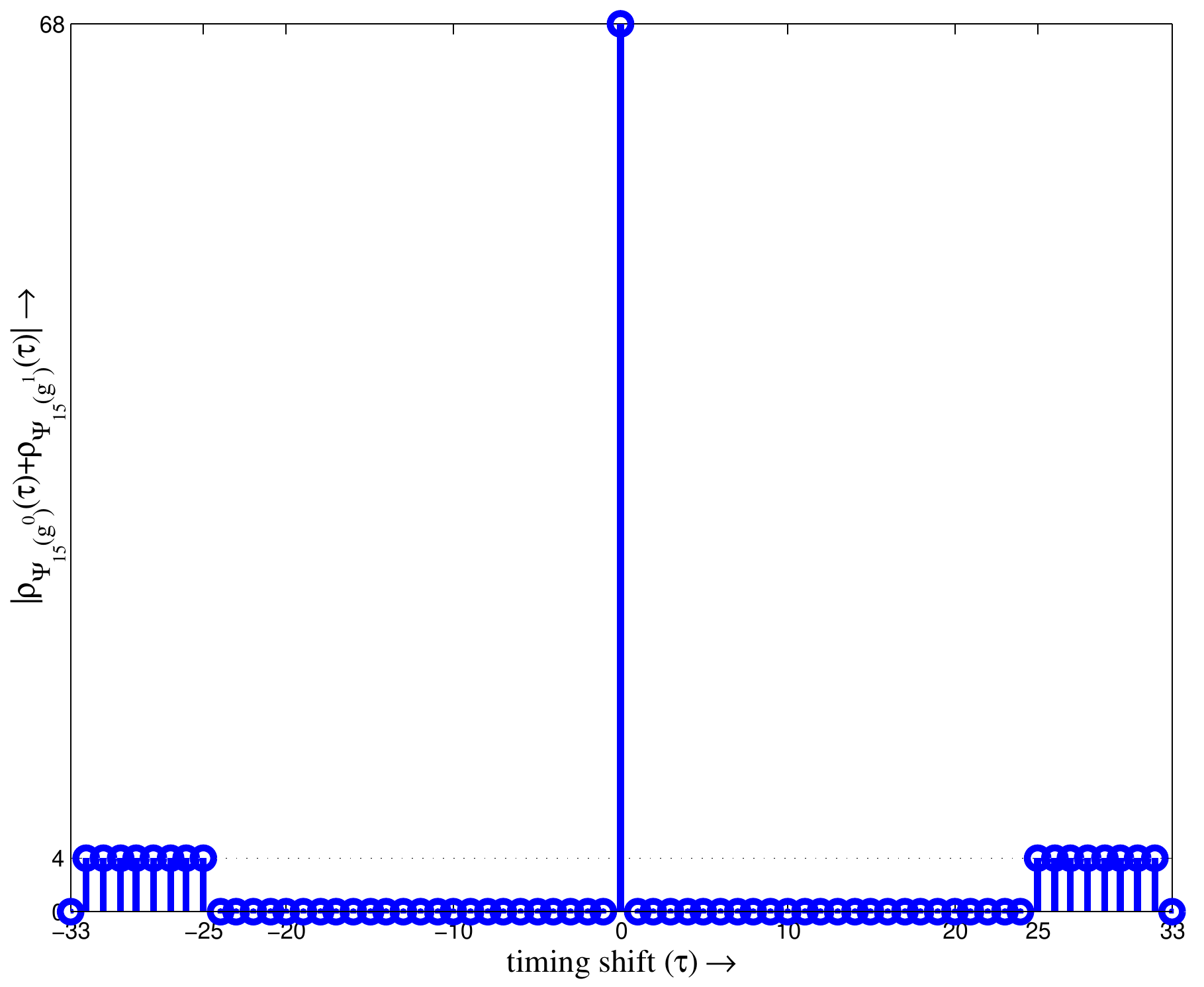}
	\caption{Autocorrelation plot of EB-ZCP in \textit{Example \ref{ex2}}.}\label{label-a}
\end{figure}

\begin{example}
	Let us consider $m=5$ and $\{\pi(0),\pi(1),\pi(2)\}=\{1,2,0\}$. For $d\in \mathbb{Z}_2$, let the Boolean function $g^d$ be given by 
	\begin{equation}
g^d=x_3\bar{x}_4\cdot \zeta^d+\bar{x}_3x_4\cdot \eta^d+d\bar{x}_3\bar{x}_4+x_3x_4,
	\end{equation}
	where
	\begin{equation}
		\zeta^d=x_1x_2+x_2x_0+dx_0,
	\end{equation}
	and
	\begin{equation}
	\eta^d=\bar{x}_1\bar{x}_2+\bar{x}_2\bar{x}_0+\bar{d}x_0.
	\end{equation}
	Since $\pi(2)=0$, the pair $(\mathbf{a},\mathbf{b})\equiv \left(\Psi_{7}(g^0),\Psi_{7}(g^1)\right)$ gives an $18$ length EB-ZCP, with ZCZ width $10$.
\end{example}


\subsection{Comparison with the previous works}

In Table \ref{table1}, we compare our proposed construction with the existing works, where the EB-ZCPs have a length of the form $2^{m-1}+2$. As we can see, our proposed construction is direct, as it does not require any sequences  at the initial stage of construction. And the resultant EB-ZCPs have the ZCZ ratio of $3/4$. 

\begin{table}[!h]
	\small
	\renewcommand{\arraystretch}{1.3}
	\resizebox{\textwidth}{!}{
		\begin{tabular}{|l|l|l|l|}
			\hline
			\begin{tabular}[c]{@{}l@{}}Construction\\ of EB-ZCPs\end{tabular}                          & Method  &\begin{tabular}[c]{@{}l@{}}Direct/ \\ Indirect \\ Construction\end{tabular}                                                                                                                       & \begin{tabular}[c]{@{}l@{}}ZCZ \\ Ratio\end{tabular} \\ \hline
			\cite{zilong_ebzcp}      & \begin{tabular}[c]{@{}l@{}}Based on truncation\\ of certain GCPs\end{tabular} & Indirect                                                       & $\frac{2}{3}$     \\ \hline 
			\cite{chen}     & Based on GBFs & Direct                                                       & $\frac{2}{3}$     \\ \hline
			\cite{Avik_ebzcp}          & \begin{tabular}[c]{@{}l@{}} Applying Insertion\\  method on certain\\GCPs\end{tabular} & Indirect & $\frac{3}{4}$     \\ \hline 
			Proposed & Based on GBFs & Direct                                       & $\frac{3}{4} $     \\ \hline
	\end{tabular}}
	\caption{Comparison with the existing works, for EB-ZCPs of length of the form $2^{m-1}+2$. \label{table1} }
\end{table}

\section{Conclusion}
In this work, a direct construction of EB-ZCPs is proposed using GBFs. The proposed EB-ZCPs are of lengths $2^{m-1}+2$ with flexible ZCZ widths of $2^{m-2}+2^{\pi(m-3)}+1$. When $\pi(m-3)=m-3$, the ZCZ ratio of the proposed EB-ZCPs are approximately equal to $3/4$, which is larger than the ZCZ ratio of the EB-ZCPs proposed by Chen.


\begin{thebibliography}{20}
	\bibitem{golay1}
	M. J. E. Golay, ``Static multislit spectrometry and its application to the
	panoramic display of infrared spectra," \emph{J. Opt. Soc. Am.}, vol. 41, no. 7,
	pp. 468-472, Jul. 1951.


\bibitem{popovic}
B. M. Popovic, ``Synthesis of power efficient multitone signals with flat
amplitude spectrum," \emph{IEEE Trans. Commun.}, vol. 39, no. 7, pp. 1031-
1033, Jul. 1991.

\bibitem{davis}
J. A. Davis and J. Jedwab, ``Peak-to-mean power control in OFDM,
Golay complementary sequences and Reed-Muller codes," \emph{IEEE Trans.
	Inf. Theory}, vol. 45, pp. 2397-2417, Nov. 1999.

\bibitem{paterson}
K. G. Paterson, ``Generalized Reed-Muller codes and power control
in OFDM modulation," \emph{IEEE Trans. Inf. Theory}, vol. 46, no. 1,
pp. 104-120, Jan. 2000.

\bibitem{doppler_gcp}
A. Pezeshki, A. R. Calderbank, W. Moran, and S. D. Howard, ``Doppler
resilient Golay complementary waveforms," \emph{IEEE Trans. Inf. Theory},
vol. 54, no. 9, pp. 4254-4266, Sep. 2008.


\bibitem{spasojevic}
P. Spasojevic and C. N. Georghiades, ``Complementary sequences for
ISI channel estimation," \emph{IEEE Trans. Inf. Theory}, vol. 47, no. 3, pp.
1145-1152, Mar. 2001.


\bibitem{fan}
P. Fan, W. Yuan, and Y. Tu, ``Z-complementary binary sequences," \emph{IEEE
	Signal Process. Lett.}, vol. 14, no. 8, pp. 509-512, Aug. 2007.

\bibitem{li2010}
X. Li, P. Fan, X. Tang, and Y. Tu, ``Existence of binary Z-complementary pairs," \emph{IEEE Signal Process. Lett.}, vol. 18, no. 1, pp.
6366, Jan. 2011.

\bibitem{zilong_obzcp}
Z. Liu, U. Parampalli, and Y. L. Guan, ``Optimal odd-length binary
Z-complementary pairs," \emph{IEEE Trans. Inf. Theory}, vol. 60, no. 9, pp.
5768-5781, Sep. 2014.

\bibitem{Avik_iwsda}
A. R. Adhikary, S. Majhi, Z. Liu, and Y. L. Guan, ``New optimal binary Z-complementary pairs of odd lengths," \emph{in Proc. The 8th IEEE International Workshop on Signal Design and its Applications in Communications}, Sapporo, Japan, pp. 14-18, Sept. 2017.

\bibitem{Avik_cs}
A. R. Adhikary and S. Majhi, ``New constructions of complementary sets of sequences of lengths non-power-of-two," in \emph{IEEE Commun. Lett.}, Early access.

\bibitem{shibu1}
S. Das, S. Majhi and Z. Liu, ``A novel class of complete complementary codes and their applications for APU matrices," in \emph{IEEE Signal Process. Lett.,} vol. 25, no. 9, pp. 1300-1304, Sept. 2018.

\bibitem{shibu2}
S. Das, S. Majhi, S. Budišin and Z. Liu, ``A new construction framework for polyphase complete complementary codes with various lengths," in \emph{IEEE Trans. Signal Process.,} vol. 67, no. 10, pp. 2639-2648, 15 May15, 2019.

\bibitem{palash1}
P. Sarkar, S. Majhi and Z. Liu, ``Optimal $Z$ -complementary code set from generalized Reed-Muller codes," in \emph{IEEE Trans. Commun.,} vol. 67, no. 3, pp. 1783-1796, March 2019.

\bibitem{palash2}
P. Sarkar, S. Majhi, H. Vettikalladi and A. S. Mahajumi, ``A direct construction of inter-group complementary code set," in \emph{IEEE Access,} vol. 6, pp. 42047-42056, 2018.


\bibitem{zilong_ebzcp}
Z. Liu, U. Parampalli, and Y. L. Guan, ``On even-period binary Z-complementary pairs with large ZCZs," \emph{IEEE Signal Process. Lett.}, vol. 21, pp.
284-287, Jun. 2014.


\bibitem{chen}
Chao-Yu Chen, ``A novel construction of Z-complementary pairs
based on generalized Boolean functions," \emph{IEEE Signal Process. Lett.}, vol. 24, pp.
284-287, Jul. 2017.








\bibitem{Avik_ebzcp}
A. R. Adhikary, S. Majhi, Z. Liu and Y. L. Guan, ``New sets of even-length binary Z-complementary pairs with asymptotic ZCZ ratio of $3/4$," \emph{IEEE Signal Process. Lett.}, vol. 25, no. 7, pp. 970-973, July 2018.


\bibitem{rathinakumar}
 A. Rathinakumar A. K. Chaturvedi, ``Complete mutually orthogonal Golay complementary sets from Reed-Muller codes" \emph{IEEE Trans. Inf. Theory} vol. 54 no. 3 pp. 1339-1346 Mar. 2008. 

\end{thebibliography}
\end{document}